\documentclass[12pt,oneside,notitlepage]{article}
\usepackage{amsfonts}
\usepackage{amsthm}
\usepackage{amsmath}
\usepackage{amscd}
\usepackage{amssymb}
\usepackage{makeidx}
\usepackage{enumerate}
\usepackage{ulem}
\usepackage{eufrak}
\usepackage{fontenc}
\usepackage[all]{xy}
\usepackage{array}
\usepackage{mathrsfs}
\usepackage{graphicx}
\usepackage[all]{xy}
\usepackage{setspace}
\usepackage{anysize}
\marginsize{1.5in}{1in}{1in}{1in}

\def\ring#1{\ifmmode \mathaccent'027 #1\else \rm\accent'027 #1\fi}

\newcommand{\RR}{\mathbb{R}}

\newcommand{\NN}{\mathbb{N}}

\newtheorem{thm}{Theorem}
\newtheorem{lemma}{Lemma}

\newtheorem{cons}{Consequence}
\newtheorem{defin}{Definition}

\newtheorem{claim}{Claim}

\begin{document}
\setcounter{page}{1}

\title{Stochastic Homology. Reduction Formulas for Computing Stochastic Betti Numbers
of Maximal Random Complexes with Discrete Probabilities.  Computation and Applications}

\date{}
\author{Todor Todorov}
\maketitle

\begin{abstract}

Given a chain complex with the only modification that each cell of the complex has a probability distribution assigned.  We will call this complex - a random complex and what should be understood in practice, is that we have a classical chain complex whose cells appear and disappear according to some probability distributions.  In this paper, we will try to find the stochastic homology of random complex, whose simplices have independent discrete distributions.

\end{abstract}

\tableofcontents

\vspace{20mm}

The development of the computer technologies nowadays is so vast that basically almost everyone in the world relies on some sort of computational device.  Devices such as personal computers, cell phones, global positioning systems, etc are part of our life, mostly because of their ability to perform millions of complex mathematical operations within a second.  Even so, the human brain is still unbeatable with its property to successfully approximate solutions of complex problems.  Without any equivocations, here is what I mean:

Everybody has seen a satellite map on Google Earth like the one of the Uptown campus of Tulane University below.

\includegraphics[scale=0.43]{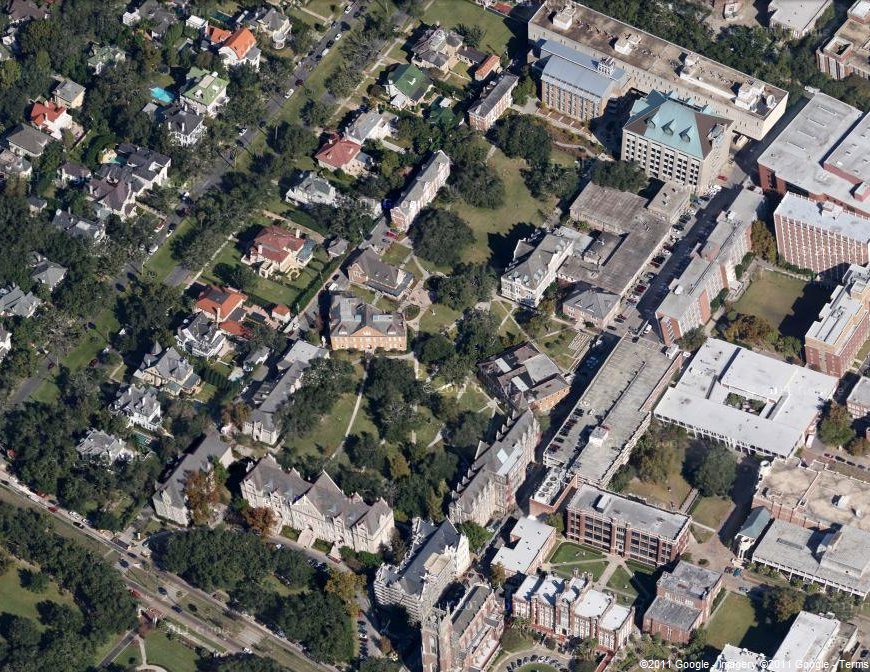}
\begin{center}
\small{Tulane's Uptown campus}
\end{center}

One very significant feature of most of the university campuses is the net of paths between the buildings like the one on Tulane's campus.  Here is another picture - this time from the area where I currently live. 
 
\includegraphics[scale=0.34]{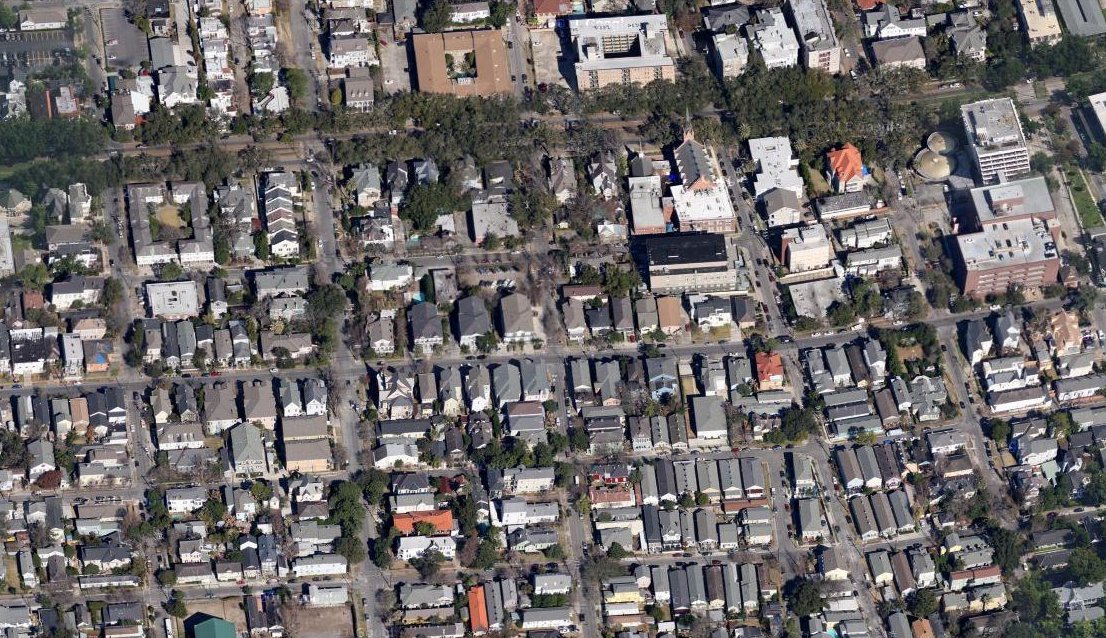}
\begin{center}
\small{Uptown New Orleans}
\end{center} 

Definitely, if one sees two pictures - one showing the net of paths on Tulane's campus and one which was taken anywhere in the residential area, he can distinguish which one is from the campus and which one is not by simply observing the existence of the gray net of paths on one of them.  Notice that we do not restrict ourselves to the angle which captures the tracks on the campus, nor the scale.  This uncountable set of pictures showing gray paths and green grass can be identified by a simple topological invariant namely

\begin{center}
\begin{tabular}{l l}
$b_0^{Green} = 22,$ & $b_1^{Green} = 0$\\
$b_0^{Gray} = 1,$ & $b_1^{Gray} = 22.$
\end{tabular}
\end{center}

Here we suppose that the number of cycles which the paths create are 22, all paths enclose grass regions and are connected.  So any picture which has more than a few gray cycles, can be considered from the Tulane's campus.  Of coarse, this might not be true at all in general, but if the maps which are covered are small enough and the pictures have enough details, there are certainly uniqueness conditions which imply that a picture belongs to a certain region.  It all seems so easy - get two satellite pictures from the regions you want to recognize, threshold them by subdividing the RGB cube into smaller cubes, get a representation of each region in terms of betti numbers corresponding to certain colors using techniques like persistence homology, do the same for the picture which was taken from the ground and hope that you are lucky enough to get unique sets of betti numbers in your picture, so that you can guess which region it corresponds to.   However, once you get your hands dirty with the real data, it turns out that it is not THAT easy.  In general, each picture contains so much noise that it is very difficult to get the actual betti numbers.  Trees, cars, people, and shadows are just the starting point of all troubles which one encounters.  Take a look at the next picture

\includegraphics[scale=0.34]{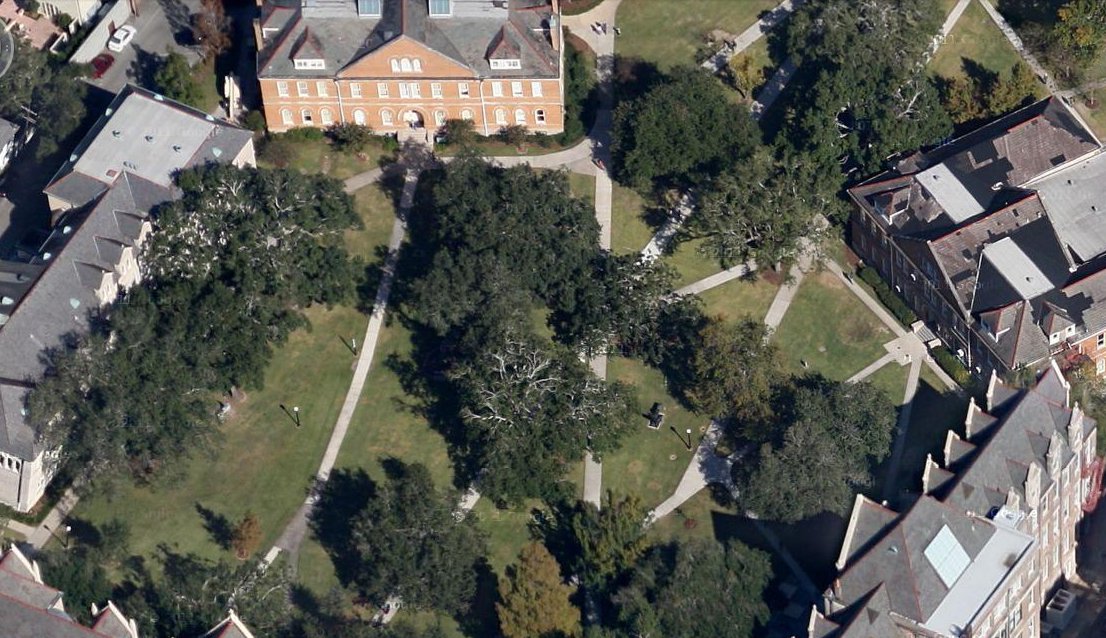}
\begin{center}
\small{Close capture of Tulane's Uptown campus}
\end{center} 
and try to count the gray cycles.  You should be able to count at least 10 of them, and all of them are approximations made by your brain.  On the other side, trying to compute the number of cycles by subdivision of the RGB cube gives terrible results no matter what subdivision of the RGB cube you take.  Thus we need a better tool to capture the topological data in the picture.  Something which would allow us to guess whether a darker region really exists, or belongs to a near lighter region, or is simply a person, etc.  We need methods which give freedom to eliminate noise, no matter how big the source is.  This brings us to the idea of stochastic homology.

Loosely speaking, stochastic homology is an extension of the idea of homology.  We build our theory like it is done in the classical homology theory - by looking at cell complexes and boundary maps.  This time, however, the fundamental unit is called random complex, which is a set of random cells.  Each random cell differs from the classical cell by the property that it has a probability distribution assigned to it.  In this paper we would only consider random complexes with discrete probabilities.  You can think of those probabilities as the probabilities of existence of each cell.  Probability of 1/2 of a certain cell means that the cell exists only during half of the time.  We would also suppose that all cells are independently distributed as dependence between the cells is the same as looking at the cells as one, and thus is not very interesting to consider.  Also, we would only consider maximal complexes, as any complex sits in a maximal.

\section{The First Example}

Suppose we are given a very simple chain complex of two points - $1$ and $2$, and an edge joining them.  As noted before, the random complex has assigned probabilities to its cells and let us denote the probabilities of the vertices $1$ and $2$ with $p_1$ and $p_2$ correspondingly, and the probability of the edge with $p_{12}$. Let's assume for completeness that $p_1=1/2$, $p_2=1/4$ and $p_{12}=1/3$.  Thus point $1$ appears only during half of the time, point $2$ appears in 1/4 of the time and the edge between them, appears in 1/3 of the time when possible.  This means that the random complex consists of two points connected with an edge in 1/24( = 1/2 * 1/3 * 1/4) of the time, there are certain moments when the complex is represented by the two points only and that happens in 1/12 ( = 1/2 * (1 - 1/3) * 1/4) of the time, in 3/8 ( = 1/2 * (1 - 1/4)) we can see the existence of point $1$, in 1/8 ( = (1 - 1/2) * 1/4) time we can observe point $2$ only and of course, in the remaining time of 3/8 ( = (1 - 1/2) * (1 - 1/4)) non of the cells exist and the complex is represented by the empty set.  

\begin{center}
\includegraphics[scale=0.50]{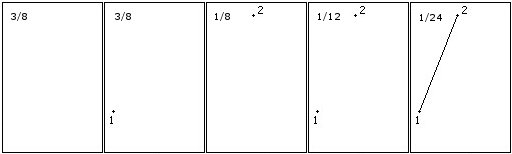}
\begin{center}
\small{Subcomplexes of the random complex}
\end{center} 
\end{center}
We already split the random complex into subcomplexes in the classical sense and to each subcomplex we can assign probability or time of existence.  Each subcomplex has certain classical topological invariants assigned like euler characteristic, homology, cohomology, homotopy groups, etc.  Informally, we can define the class of expected invariants as
$$expected \; invariant = \sum_{\Delta} value \; of \; the \; classical \; invariant(\Delta) * probability(\Delta) $$
where $\Delta$ runs over all possible chain subcomplexes. Using the last formula, it is very easy to compute the expected number of components, denote by $b_0^E$, of the above complex and it is given by
$$b_0^E = 0 * 3/8 + 1 * 3/8 + 1 * 1/8 + 2 * 1/12 + 1 * 1/24 = 17/24.$$
Of course, in this simple case the expected number of components is equal to the expected Euler characteristic, denoted as $\chi^E$, so $\chi^E = 17/24$ also.

Although already mentioned, let us give the following

\begin{defin}
The expected k-th betti number $b_k^E$ is the expectation of the classical k-th betti number over all possible configurations of complexes, i.e.
$$b_k^E = \sum_{\Delta} b_k(\Delta) p(\Delta)$$
where $\Delta$ runs over all possible chain complexes.
\end{defin}

Consider the previous example and let us try to find a formula for $b_0^E$ in terms of the probabilities $p_1, p_2$ and $p_{12}$.  Using the above arguments, we have
$$\begin{array}{rcl}
b_0^E &=& 0 * (1 - p_1) * (1 - p_2) + 1 * p_1 * (1 - p_2) + 1 * (1 - p_1) * p_2\\
& & + 2 * p_1 * p_2 * (1 - p_{12}) + 1 * p_1 * p_2 * p_{12} = p_1 + p_2 - p_1 * p_2 * p_{12}.
\end{array}$$
Note that the probability $p_{12}$ appears only when both points exist in the complex.  Also, for the sake of short notations we would only write $p_{12}$  for the term $p_1 * p_2 * p_{12}$ or it should be understood that the probabilities of the lower dimensional cells are already implemented into the probability of the higher dimensional cell.  There is a nice geometric representation for such formulas and we would prefer to utilize it whenever possible.  The coefficient of each summand appears in the upper left corner.

\begin{center}
\includegraphics[scale=0.50]{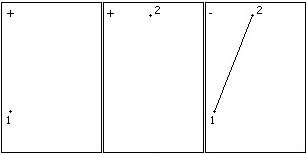}
\begin{center}
\small{Geometric representation of $b_0^E$ of a maximal random complex over two points}
\end{center} 
\end{center}
Actually, if someone simplifies the expected zeroth betti number of the maximal random complex over 3 points, it turns out that
$$b_0^E = p_1 + p_2 + p_3 - p_{12} - p_{13} - p_{23} + p_{12}*p_{13}*p_{23}$$

\begin{center}
\includegraphics[scale=0.50]{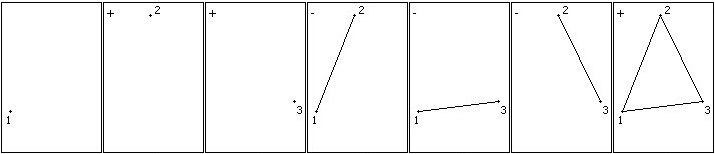}
\begin{center}
\small{Geometric representation of $b_0^E$ of a maximal random complex over three points}
\end{center} 
\end{center}
and it seems that we can already guess what the general formula for the expected zeroth betti number of the maximal random complex over n points could be.  Well, not quite!  In fact, the formula for the maximal 4 points complex is given 

\begin{center}
\includegraphics[scale=0.45]{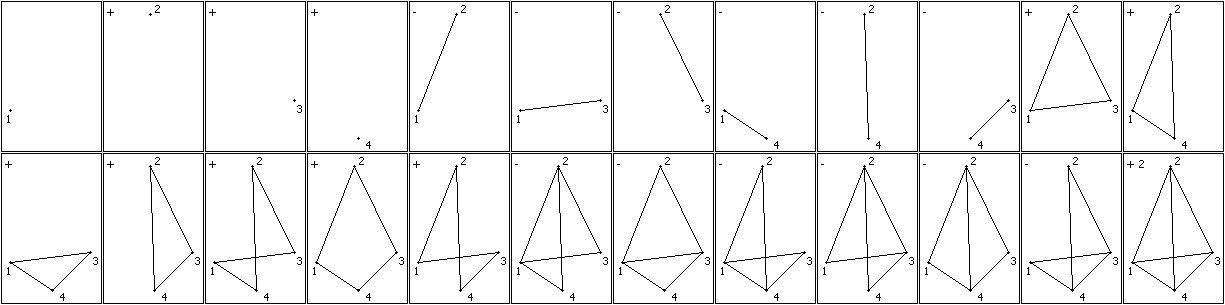}
\begin{center}
\small{Geometric representation of $b_0^E$ of a maximal random complex over four points}
\end{center} 
\end{center}
and the one for 5 points is given by

\begin{center}
\includegraphics[scale=0.45]{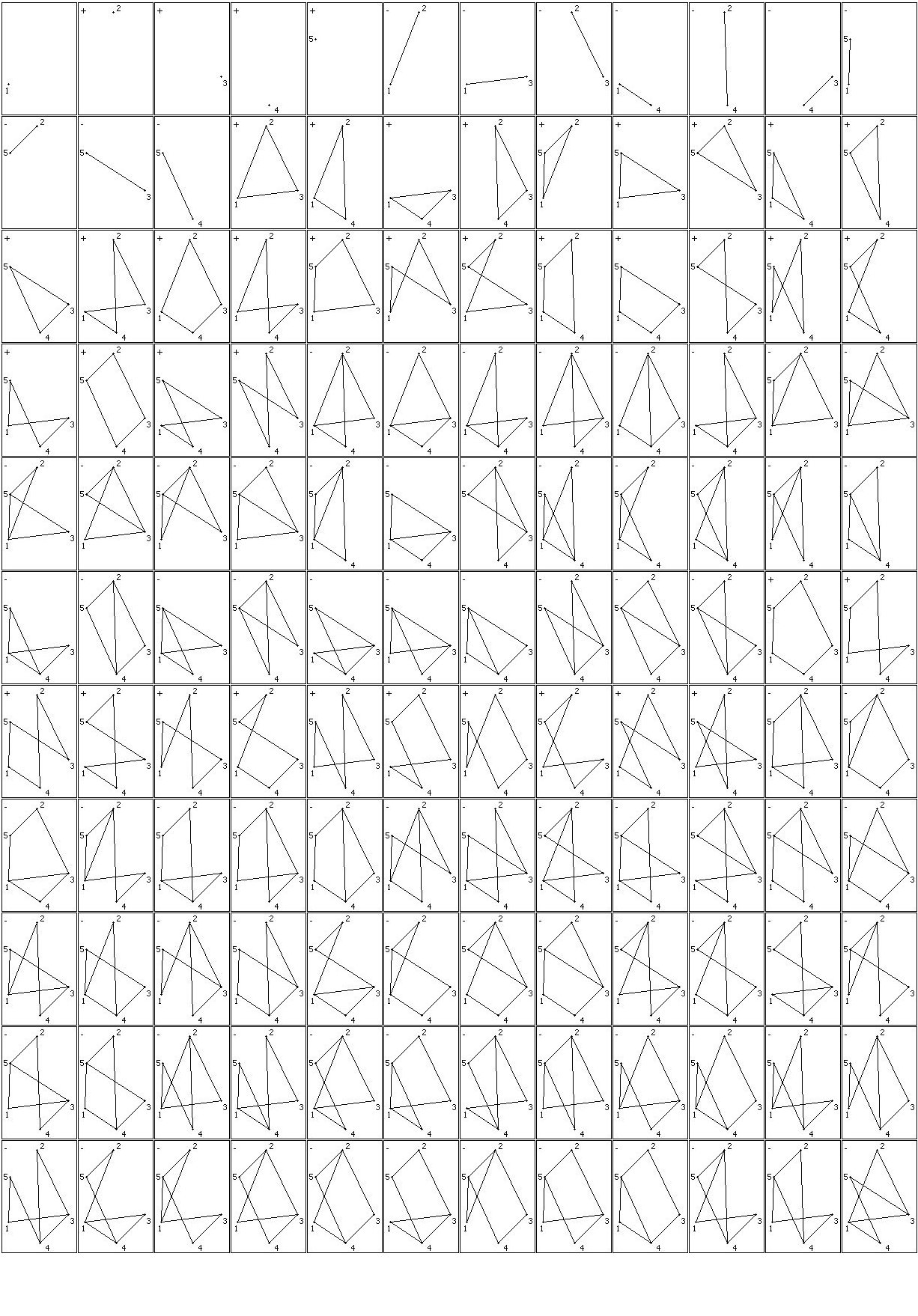} 
\end{center}

\begin{center}
\includegraphics[scale=0.45]{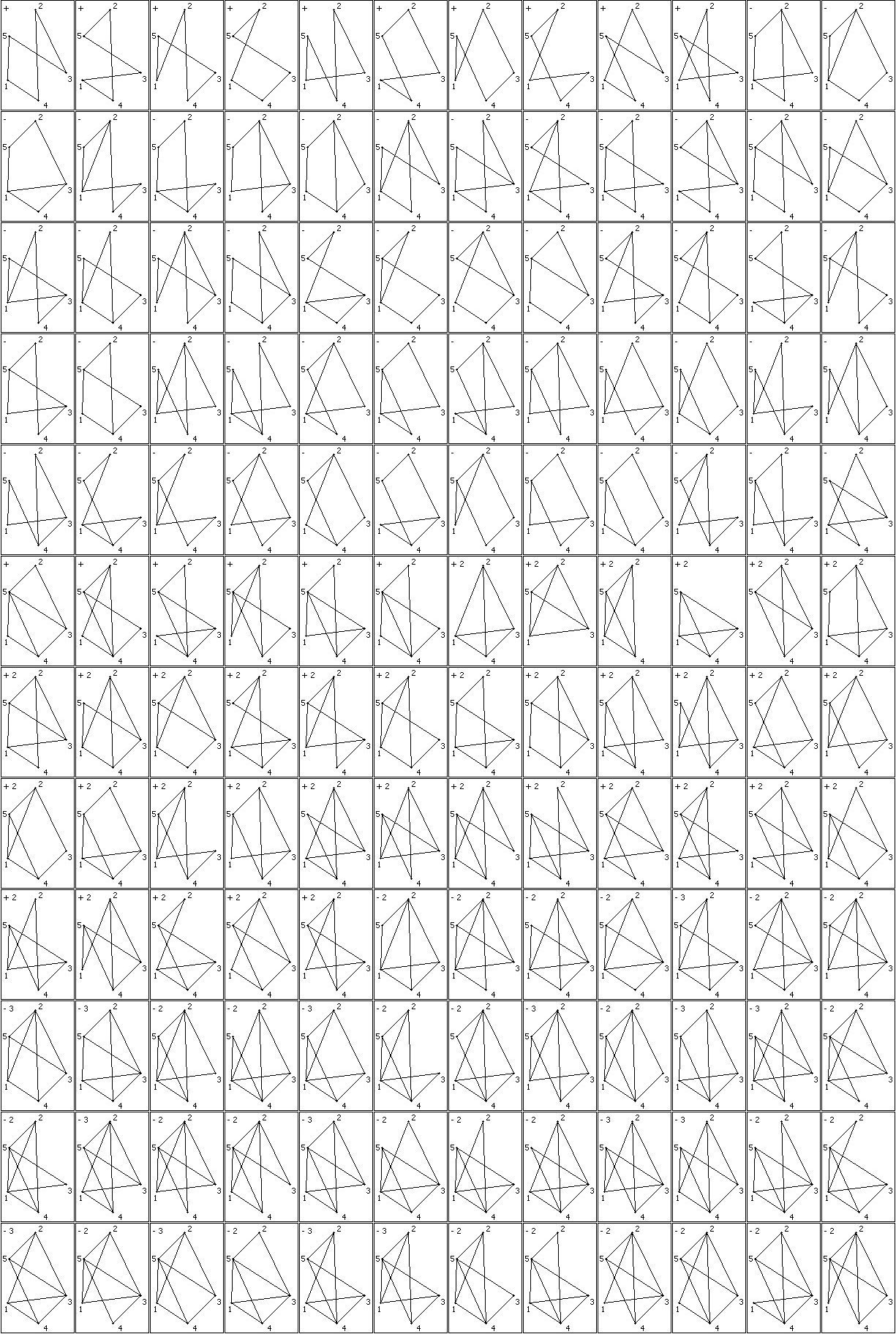}
\end{center}

\begin{center}
\includegraphics[scale=0.45]{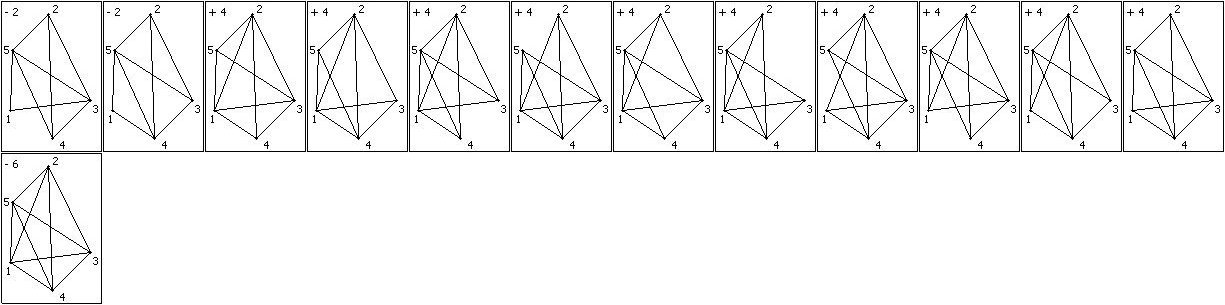}
\begin{center}
\small{Geometric representation of $b_0^E$ of a maximal random complex over four points}
\end{center}
\end{center}

The formula for the maximal random complex over six points include 12987 terms.  However, the total number of subcomplexes of the maximal complex over $n$ points is given by 
$$\sum_{k=0}^n \binom{n}{k} 2^k.$$
The following table represents the last number
$$\begin{array}{llll}
Points: & 1 & \quad Complexes: & 2\\
& 2 & & 5\\
& 3 & & 18\\
& 4 & & 113\\
& 5 & & 1450\\
& 6 & & 40069\\
& 7 & & 2350602\\
& 8 & & 286192513\\
& 9 & & 2494306930\\
\end{array}$$
Among this mess of summands, we will try to find a recurrence.

\section{Term Structure of $b_0^E$}

The first natural question which we may ask is about the expected zeroth betti number of n points and no edges in the random complex.  The answer is given by the next

\begin{claim}\label{cl:one}
Let the random complex $X$ consists of n points only.  Then zeroth betti of $X$ is   
$$b_0^E(X)= p_1 + \ldots + p_n$$
\end{claim}

\textbf{Proof.}  If we sum over the number of points, i.e. over the number of connected components, we have the following identity
\begin{equation}\label{eq:fla1}
\begin{array}{lcr}
\sum_{i=1}^n p_i (1-p_1) \ldots \widehat{(1-p_i)} \ldots (1-p_n) & & \\
+ 2 \sum_{i<j} p_i p_j (1-p_1) \ldots \widehat{(1-p_i)} \ldots \widehat{(1-p_j)} \ldots (1-p_n) & & \\
+ \ldots & & \\
+ n p_1 \ldots p_n & = & p_1 + \ldots + p_n
\end{array}
\end{equation}
where the coefficient in front of each sum is the betti number of the corresponding complexes.  The proof of equation (\ref{eq:fla1}) is done by induction.  Obviously, for $n=1$ the formula is true.  Suppose that equation (\ref{eq:fla1}) is true for $n$ and consider the following sum
$$\begin{array}{lcr}
\sum_{i=1}^{n+1} p_i (1-p_1) \ldots \widehat{(1-p_i)} \ldots (1-p_{n+1}) & & \\
+ 2 \sum_{i<j} p_i p_j (1-p_1) \ldots \widehat{(1-p_i)} \ldots \widehat{(1-p_j)} \ldots (1-p_{n+1}) & & \\
+ 3 \sum_{i<j<k} p_i p_j p_k (1-p_1) \ldots \widehat{(1-p_i)} \ldots \widehat{(1-p_j)} \ldots \widehat{(1-p_k)} \ldots (1-p_{n+1}) & & \\
+ \ldots & & \\
+ n \sum_{i=1}^{n+1} (1-p_i) p_1 \ldots \widehat{p_i} \ldots p_{n+1} & & \\
+(n+1) p_1 \ldots p_{n+1} & &\\
= (1-p_{n+1}) \sum_{i=1}^{n} p_i (1-p_1) \ldots \widehat{(1-p_i)} \ldots (1-p_n) + p_{n+1} (1-p_1) \ldots (1-p_n)& & \\
+ 2 (1-p_{n+1}) \sum_{i<j} p_i p_j (1-p_1) \ldots \widehat{(1-p_i)} \ldots \widehat{(1-p_j)} \ldots (1-p_n)& & \\
+ 2 p_{n+1} \sum_{i=1}^{n} p_i (1-p_1) \ldots \widehat{(1-p_i)} \ldots (1-p_n) & & \\
+ 3 (1-p_{n+1}) \sum_{i<j<k} p_i p_j p_k (1-p_1) \ldots \widehat{(1-p_i)} \ldots \widehat{(1-p_j)} \ldots \widehat{(1-p_k)} \ldots (1-p_{n+1})& & \\
+ 3 p_{n+1} \sum_{i<j} p_i p_j (1-p_1) \ldots \widehat{(1-p_i)} \ldots \widehat{(1-p_j)} \ldots (1-p_n) & & \\
+ \ldots & & \\
+ (1-p_{n+1}) n p_1 \ldots p_n + n p_{n+1} \sum_{i=1}^{n} (1-p_i) p_1 \ldots \widehat{p_i} \ldots p_{n+1} & & \\
+(n+1) p_1 \ldots p_{n+1} & &\\
= (1-p_{n+1})(p_1 + \ldots + p_n) + p_{n+1}(p_1 + \ldots + p_n) + p_{n+1} (1-p_1) \ldots (1-p_n)& & \\
+ p_{n+1} \sum_{i=1}^{n} p_i (1-p_1) \ldots \widehat{(1-p_i)} \ldots (1-p_n) & & \\
+ p_{n+1} \sum_{i<j} p_i p_j (1-p_1) \ldots \widehat{(1-p_i)} \ldots \widehat{(1-p_j)} \ldots (1-p_n) & & \\
+ \ldots & & \\
+ p_{n+1} \sum_{i=1}^{n} (1-p_i) p_1 \ldots \widehat{p_i} \ldots p_{n+1} & & \\
+ p_1 \ldots p_{n+1} & &\\
= (p_1 + \ldots + p_n) + p_{n+1} (1-p_1) \ldots (1-p_n)& & \\
+ p_{n+1} \sum_{i=1}^{n} p_i (1-p_1) \ldots \widehat{(1-p_i)} \ldots (1-p_n) & & \\
+ p_{n+1} \sum_{i<j} p_i p_j (1-p_1) \ldots \widehat{(1-p_i)} \ldots \widehat{(1-p_j)} \ldots (1-p_n) & & \\
+ \ldots & & \\
+ p_{n+1} \sum_{i=1}^{n} (1-p_i) p_1 \ldots \widehat{p_i} \ldots p_{n+1} & & \\
+ p_{n+1} p_1 \ldots p_n & &\\
\end{array}$$

In the last equation we split the sum into two sums and use the induction hypothesis.  The Claim follows by the next result.

\begin{flushright}
$\Box$
\end{flushright}

\begin{claim}  The following equation holds
$$\begin{array}{lcr}
(1-p_1) \ldots (1-p_n)+ \sum_{i=1}^{n} p_i (1-p_1) \ldots \widehat{(1-p_i)} \ldots (1-p_n)& & \\
 + \sum_{i<j} p_i p_j (1-p_1) \ldots \widehat{(1-p_i)} \ldots \widehat{(1-p_j)} \ldots (1-p_n) + \ldots + p_1 \ldots p_n = 1
\end{array}$$
\end{claim}

\textbf{Proof.} Using induction again, we have that the result holds for $n=1$.  Suppose it is true for $n$ and split the sum as before
$$\begin{array}{lcr}
(1-p_1) \ldots (1-p_{n+1})+ \sum_{i=1}^{n+1} p_i (1-p_1) \ldots \widehat{(1-p_i)} \ldots (1-p_{n+1})& & \\
 + \sum_{i<j} p_i p_j (1-p_1) \ldots \widehat{(1-p_i)} \ldots \widehat{(1-p_j)} \ldots (1-p_{n+1}) + \ldots + p_1 \ldots p_{n+1}& & \\
= (1-p_{n+1}) (1-p_1) \ldots (1-p_n) + (1-p_{n+1}) \sum_{i=1}^{n} p_i (1-p_1) \ldots \widehat{(1-p_i)} \ldots (1-p_n)& & \\
 + p_{n+1} (1-p_1) \ldots (1-p_n) + (1-p_{n+1}) \sum_{i<j} p_i p_j (1-p_1) \ldots \widehat{(1-p_i)} \ldots \widehat{(1-p_j)} \ldots (1-p_n) + \ldots& & \\
+ p_1 \ldots p_{n+1} = (1-p_{n+1}) + p_{n+1} = 1
\end{array}$$
\begin{flushright}
$\Box$
\end{flushright}

Actually, this long proof could be omitted if we notice that since there are no edges in the randoms complex consisting of n points, then these points are independent variables and thus
$$E(\sum_{i=1}^n p_i) = \sum_{i=1}^n E(p_i)$$
which is just another way of writing Claim \ref{cl:one}.

Next, we concentrate on a more general random complex $X$ over n points which has cells in higher dimensions.  The first observation is that $b_0^E$ depends only on the zeroth and first chain subcomplexes of $X$ as the classical betti number depends only on those subcomplexes.  Then while we sum all various monomials corresponding to the subcomplexes of $X$, we can note that basically there are two types of monomials - one corresponding to complexes in which there exist no edges
$$p_{i_1} \ldots p_{i_p} (1 - p_{j_1 j_2}) \ldots (1 - p_{j_q j_r})$$
and another in which there is at least one edge
$$p_{i_1} \ldots p_{i_p} (1 - p_{j_1 j_2}) \ldots (1 - p_{j_q j_r}) p_{k_1 k_2} \ldots p_{k_s k_t}.$$
If we consider all summands from the first type and expand all the products in the parentheses, there will be summands of type
$$p_{i_1} \ldots p_{i_p}.$$
As we know from Claim \ref{cl:one}, the sum of all these summands is $p_1 + \ldots + p_n$, thus we can split the expected zeroth betti of $X$ as
$$b_0^E (X) = p_1 + \ldots + p_n + \mbox{(terms having a variable corresponding to an edge)}.$$

Now we put our efforts on deciphering the second part of the polynomial $b_0^E(X)$.  We try to figure out the coefficient in front of a general monomial containing an edge in the next

\begin{lemma}
Let $m$ be the minimum amount of points which are covered by the edges $p_{i_1 i_2}, \ldots, p_{i_{2n-1} i_{2n}}$ and $p_1, \ldots, p_m$ are the probabilities of the covered vertices.  Then the coefficient in front of $p_1 \ldots p_m p_{m+1} \ldots p_k p_{i_1 i_2} \ldots p_{i_{2n-1} i_{2n}}$ for $n>0$ in the reduced polynomial $b_0^E$ is given by 
\begin{equation}\label{eq:equ2}
\sum_{i=0}^{k-m} \binom{k-m}{i} (-1)^i \sum_{\substack{| \omega | = 0 \\ \omega \in C(i_1, \ldots, i_{2n})}}^n (-1)^{| \omega |}(k-i-n+|\omega|+b_1^\omega)
\end{equation}
Here $\omega$ is a combinatorial term which runs over the edges $p_{i_1 i_2} \ldots p_{i_{2n-1} i_{2n}}$ and builds different subcomplexes $X^\omega$.  $|\omega|$ is the number of non-existing edges and $b_1^\omega$ is the first betti number corresponding to the complex $X^\omega$.
\end{lemma}

\textbf{Proof.}  Since the number of vertices $m$ is the minimum number of points which are covered by the edges then $p_1 \ldots p_m$ should always exist in all subcomplexes as if one of them fails to exist, then there will be no corresponding edge term in the monomial.  In this sense, $p_1 \ldots p_m$ are fixed to exist.  While, the terms $p_{m+1} \ldots p_k$ have no influence on the edges, we can freely choose the points $m+1, \ldots, k$ to exist or not.  Same is true for the edges corresponding to $p_{i_1 i_2} \ldots p_{i_{2n-1} i_{2n}}$.  There might be more edge terms, but we must choose them to fail to exist as otherwise they will contribute an unnecessary term.  Thus we need to sum all possible subcomplexes in which the points corresponding to $p_1 \ldots p_m$ exist, the points corresponding to $p_{m+1} \ldots p_k$ and the edges corresponding to $p_{i_1 i_2} \ldots p_{i_{2n-1} i_{2n}}$ may or may not exist, and all other possible edges must not exist.  The term $\binom{k-m}{i}$ in the above formula is due to the number of choices for the free points. $(-1)^i$ gives the sign after expanding, where $i$ is the number of non existing free points.  $(k-i-n+|\omega|+b_1^\omega)$ gives the zeroth betti number of a fixed subcomplex as we have $k-i$ total points, subtract the number of edges $n+|\omega|$ and add the first betti number (Euler-Poincare formula).  $(-1)^{| \omega |}$ gives the sign of non-existing edges and we have to sum over all possible existing or non-existing edges and so we get the coefficient (\ref{eq:equ2}).
\begin{flushright}
$\Box$
\end{flushright}

Our next goal is to simplify the coefficient (\ref{eq:equ2}).   First break the second sum into two sums
$$\begin{array}{lcr}
\sum_{i=0}^{k-m} \binom{k-m}{i} (-1)^i \left( \sum_{\substack{| \omega | = 0 \\ \omega \in C(i_1, \ldots, i_{2n})}}^n (-1)^{| \omega |}(k-i) + \sum_{\substack{| \omega | = 0 \\ \omega \in C(i_1, \ldots, i_{2n})}}^n (-1)^{| \omega |}(-n+|\omega|+b_1^\omega) \right)& &
\end{array}$$
The first summand does not depend combinatorially on $\omega$, so we can replace the combinatorial term by a binomial coefficient
$$\begin{array}{lcr}
\sum_{i=0}^{k-m} \binom{k-m}{i} (-1)^i \left( (k-i) \sum_{| \omega | = 0 }^n (-1)^{| \omega |}\binom{n}{| \omega |} + \sum_{\substack{| \omega | = 0 \\ \omega \in C(i_1, \ldots, i_{2n})}}^n (-1)^{| \omega |}(-n+|\omega|+b_1^\omega) \right)& & \\
= \sum_{i=0}^{k-m} \binom{k-m}{i} (-1)^i \left( (k-i)(1-1)^n + \sum_{\substack{| \omega | = 0 \\ \omega \in C(i_1, \ldots, i_{2n})}}^n (-1)^{| \omega |}(-n+|\omega|+b_1^\omega) \right)& & \\
= \sum_{i=0}^{k-m} \binom{k-m}{i} (-1)^i \sum_{\substack{| \omega | = 0 \\ \omega \in C(i_1, \ldots, i_{2n})}}^n (-1)^{| \omega |}(-n+|\omega|+b_1^\omega)& & 
\end{array}$$
The last equality follows by the fact that n is nonzero.  The sums split variables, so we can evaluate the first one separately
\begin{equation}\label{eq:equ3}
\begin{array}{lcr}
= (1-1)^{k-m} \sum_{\substack{| \omega | = 0 \\ \omega \in C(i_1, \ldots, i_{2n})}}^n (-1)^{| \omega |}(-n+|\omega|+b_1^\omega)& & \\
= \left\{ \begin{array}{lcr} \sum_{\substack{| \omega | = 0 \\ \omega \in C(i_1, \ldots, i_{2n})}}^n (-1)^{| \omega |}(-n+|\omega|+b_1^\omega)& & \mbox{if } k=m\\
0 & & \mbox{if } k \neq m
\end{array}
\right.
\end{array}
\end{equation}

If $b_1^\omega = 0$ for all $\omega$, the coefficient is equal to
$$\begin{array}{rcl}
\sum_{| \omega | = 0 }^n (-1)^{| \omega |} (-n+|\omega|) \binom{n}{|\omega|} & = & \sum_{| \omega | = 0 }^n (-1)^{| \omega |+1} n \binom{n-1}{|\omega|}\\
 & = & (-1)^n n \sum_{| \omega | = 0 }^n (-1)^{n-| \omega |+1} \binom{n-1}{|\omega|}\\
 & = & (-1)^n n (1-1)^{n-1}\\
 & = & \left\{ \begin{array}{rcr} -1 & & \mbox{if } n=1\\
0 & & \mbox{if } n \neq 1
\end{array}
\right.
\end{array}$$
Thus we proved the following

\begin{claim}
All edges come with coefficient -1.
\end{claim}

Finally, let us examine the case when $b_1 \neq 0$ in equation (\ref{eq:equ3}).  We can split the last sum ones again and substitute the combinatorial term in the second one
$$\begin{array}{lcr}
= \sum_{\substack{| \omega | = 0 \\ \omega \in C(i_1, \ldots, i_{2n})}}^n (-1)^{| \omega |}b_1^\omega - \sum_{\substack{| \omega | = 0 \\ \omega \in C(i_1, \ldots, i_{2n})}}^n (-1)^{| \omega |}(-n+|\omega|)& & \\
= \sum_{\substack{| \omega | = 0 \\ \omega \in C(i_1, \ldots, i_{2n})}}^n (-1)^{| \omega |}b_1^\omega - \sum_{| \omega | = 0}^n (-1)^{| \omega |}(-n+|\omega|) \binom{n}{|\omega|}& & \\
\end{array}$$
Now the second sum is 0 since one edge can never build a cycle, thus we have the following

\begin{thm}\label{thm:main}
Suppose that $b_1$ of the corresponding to $p_1 \ldots p_m p_{i_1 i_2} \ldots p_{i_{2n-1} i_{2n}}$ complex is nonzero. Then the coefficient in front of $p_1 \ldots p_m p_{i_1 i_2} \ldots p_{i_{2n-1} i_{2n}}$ in the reduced polynomial $b_0^E(X)$ is given by
\begin{equation}\label{eq:theq}
\sum_{\substack{| \omega | = 0 \\ \omega \in C(i_1, \ldots, i_{2n})}}^n (-1)^{| \omega |} b_1^\omega.
\end{equation}
\end{thm}
\textbf{Example.}  Set all $p_i$ and $p_{ij}$ to be either 0 or 1 and plug into the formula for $b_0^E$.  Then you get another decomposition of the classical $b_0$.

\begin{defin}
We call a simple 1-cycle a complex in which each vertex is covered by exactly two edges.  In general, simple n-cycle is a complex in which each (n-1)-cell is covered by exactly two n-cells.
\end{defin}

\begin{cons} The coefficient of the monomial corresponding to a simple 1-cycle is always 1.
\end{cons}

\textbf{Proof.}  The only complex with nonzero first betti number is the complex consisting of the simple cycle.

\begin{flushright}
$\Box$
\end{flushright}

\begin{defin}
We call a 1-spike the connected union of edges which do not constitute a 1-cycle in a complex.  Similarly, n-spike is the connected union of n-cells which do not constitute a n-cycle.
\end{defin}

Suppose now we are given a simple 1-cycle of any length with $k$ 1-spikes either attached to it or not, $k>0$.  The only way to decompose the complex so that each new complex has nonzero first betti number is to combinatorially remove all the spikes.  Using Theorem \ref{thm:main} we get that the coefficient of the monomial corresponding to the initial complex is zero since k is nonzero and
$$\sum_{i=0}^k (-1)^i \binom{k}{i} = (1-1)^k = 0.$$

\begin{defin}
Two complexes are called to be n-non-intersecting if they do not intersect in a n-cell but they may intersect in cells of lower than n dimension.  An intersection of two complexes is said to be n-nonempty if the complexes intersect in a n-cell.
\end{defin}

It turns out that there are no terms in the $b_0^E$ polynomial which correspond to two 1-non-intersecting cycles of lengths $n_1>0$ and $n_2>0$ since in this case the coefficient is given by
$$\begin{array}{lcr}\left(2 + \sum_{i=1}^{n_1} (-1)^i \binom{n_1}{i} + \sum_{i=1}^{n_2} (-1)^i \binom{n_2}{i} \right) = \left( \sum_{i=0}^{n_1} (-1)^i \binom{n_1}{i} + \sum_{i=0}^{n_2} (-1)^i \binom{n_2}{i} \right) & &\\
= \left((1-1)^{n_1}+(1-1)^{n_2}\right) = 0 
\end{array}$$

Suppose we are given a complex consisting of two cycles of lengths $n_1>0$ and $n_2>0$ with 1-nonempty intersection which occurs in $k>0$ edges and no spikes.  Then the coefficient of the corresponding monomial is -1 since
$$\begin{array}{lcr} \left(2 + \sum_{i=1}^{n_1} (-1)^i \binom{n_1}{i} + \sum_{i=1}^{n_2} (-1)^i \binom{n_2}{i} + \sum_{i=1}^{k} (-1)^i \binom{k}{i} \right)& &\\
= \left( \sum_{i=0}^{n_1} (-1)^i \binom{n_1}{i} + \sum_{i=0}^{n_2} (-1)^i \binom{n_2}{i} + \sum_{i=0}^{k} (-1)^i \binom{k}{i} -1 \right) & &\\
= \left((1-1)^{n_1}+(1-1)^{n_2}+(1-1)^k -1\right) = -1. 
\end{array}$$

\textbf{Example.}  So far we already discovered most of the formula for the expected zeroth betti number of the tetrahedron.  The only unknown which remains is the coefficient of the maximal complex which we denote with $C$.
\begin{equation}\label{eq:pyr}
\begin{array}{rcl}
b_0^E(X) & = & p_1 + p_2 + p_3 + p_4\\
 & & - p_{12} - p_{13} - p_{14} - p_{23} - p_{24} - p_{34}\\
 & & + p_{12} p_{13} p_{23} + p_{12} p_{14} p_{24} + p_{13} p_{14} p_{34} + p_{23} p_{24} p_{34}\\
 & & + p_{12} p_{23} p_{34} p_{14} + p_{12} p_{24} p_{34} p_{13} + p_{13} p_{23} p_{24} p_{14}\\
 & & - p_{12} p_{13} p_{14} p_{23} p_{24} - p_{12} p_{13} p_{14} p_{23} p_{34} - p_{12} p_{13} p_{14} p_{24} p_{34}\\
 & & - p_{12} p_{13} p_{23} p_{24} p_{34} - p_{12} p_{14} p_{23} p_{24} p_{34} - p_{13} p_{14} p_{23} p_{24} p_{34}\\
 & & + C p_{12} p_{13} p_{14} p_{23} p_{24} p_{34} 
\end{array}
\end{equation}

One way to get the coefficient $C$ is by using formula (\ref{eq:theq}) and combinatorially decomposing the corresponding complex.  We leave this as an exercise.  Another more efficient way to calculate the coefficient $C$ is to set all probabilities of the points and edges to be equal to one and substitute them in equation (\ref{eq:pyr}).  Thus $b_0^E(X)$ in this case is equal to $b_0(X) = 1$ and we can solve the equation 
$$1 = 4 - 6 + 4 + 3 - 6 + C$$
for $C$ and get $C=2$.

As one can notice, the computing coefficients is via formula (\ref{eq:theq}) is slow and resource demanding.  We will be looking for a more efficient formula for computing the coefficients of the monomials in the polynomial $b_0^E$.  For simplicity, we introduce the following notations.

The initial complex will usually be denoted by $\Delta$ and the set of all complexes which are derived from $\Delta$ by removing exactly $i$ edges will be denoted by $\Delta^i$.  Notice that $\Delta^0$ is just $\Delta$.  The notation $b_1(\Delta^i)$ will stand for the sum of the betti numbers of all complexes which belong to $\Delta^i$ which we can write as
$$b_1(\Delta^i) = \sum_{\substack{\omega \in C(i_{h_1}, \ldots, i_{h_i}) \\ | \omega | = i}} b_1^\omega.$$
Similarly, $c_1(\Delta)$ will mean the coefficient of $\Delta$ given by formula (\ref{eq:theq}) and $c_1(\Delta^i)$ will denote the sum of the coefficients of the elements of $\Delta^i$.  Here the subindex in $c_1$ is used to emphasize that the coefficient corresponds to the 1-cell complex.  Considering the results above, we artificially define $c_1(\Delta)$ to be zero if $\Delta$ has $b_1=0$.

Suppose the random complex $\Delta$ has $n$ edges.  Using the notations we just introduced, we can rewrite formula (\ref{eq:theq}) as
\begin{equation}\label{eq:sme}
c_1(\Delta) = \sum_{i=0}^n (-1)^i b_1(\Delta^i).
\end{equation}
Our next goal is to keep $b_1(\Delta^0)$ and write $b_1(\Delta^1)$ in terms of $c_1(\Delta^1)$.

The careful reader can notice that after removing an edge from the terms of $\Delta^1$, then the set of all new complexes covers $\Delta^2$ twice.  In case we remove two edges from the complexes of $\Delta^1$, then we can notice that $\Delta^3$ is covered three times and we can conjecture that for any $i$, $\Delta^i$ is covered $i$ times.  Notice that $\Delta^1$ has $\binom{n}{1}$ elements each having $n-1$ edges and $\Delta^i$ has $\binom{n}{i}$ elements each with $n-i$ edges.  There are $\binom{n-1}{1}$ choices to remove an edge from a fixed element of $\Delta^1$ and there are $n$ such elements, so $\Delta^2$ is covered twice.  Next, we have $2*\binom{n}{2}$ elements and we remove one more edge and it is easy to see that $\Delta^3$ is covered three times.  The statement is proved by induction and will be left as an exercise.

Using the last observation, we can rewrite equation (\ref{eq:sme}) as
$$\begin{array}{rcl}
c_1(\Delta) & = & b_1(\Delta^0) - b_1(\Delta^1) + \ldots + (-1)^n n b_1(\Delta^n) + (- b_1(\Delta^2) + 2 b_1(\Delta^3) - \ldots\\
& & + (-1)^{n-1} (n-1) b_1(\Delta^n))\\
& = & b_1(\Delta^0) - c_1(\Delta^1) - b_1(\Delta^2) + 2 b_1(\Delta^3) - \ldots + (-1)^{n-1} (n-1) b_1(\Delta^n).
\end{array}$$

The strategy remains the same - change the next summand $b_1(\Delta^2)$ to $c_1(\Delta^2)$.  This time we need to calculate how many times $\Delta^2$ covers $\Delta^3$.  So in order to save time and efforts, let us do it in general.  We need to figure out how many times the elements from $\Delta^j$ cover the elements in $\Delta^i$, for $j<i$.  The number of edges of each element of $\Delta^j$ is $n-j$ and we need to remove $i-j$ 1-cells in order to drop to level $\Delta^i$.  There are $\binom{n}{j}$ elements in $\Delta^j$, each covering $\binom{n-j}{i-j}$ elements from $\Delta^i$, i.e. $\binom{n}{j} \binom{n-j}{i-j}$ compared to $\binom{n}{i}$ gives that $\Delta^j$ covers $\Delta^i$ exactly $\binom{i}{j}$ times.

Thus, the coefficient of $b_1(\Delta^m)$, after changing the first $m-1$ betti terms to the corresponding coefficients $c_1$, is always
$$(-1)^m \left( 1 - \binom{m}{1} + \binom{m}{2} - \binom{m}{3} + \ldots + (-1)^{m-1} \binom{m}{m-1} \right) = (-1)^m (-1)^{m+1} = -1.$$

We just proved the next
\begin{thm}\label{thm:govna}
Using the above notations and definitions, the following formula holds
\begin{equation}\label{eq:qko}
\sum_{i=0}^n c_1(\Delta^i) = b_1(\Delta).
\end{equation}
\end{thm}
It is a good moment to verify this formula by an

\textbf{Example.} $c_1(Pyramid) = 3 - 6(-1) - 3(1) - 4(1) = 2$.

So far, we know that any 1-cycle with $k$ 1-spikes have a coefficient equal to 0.  We also know that the coefficient of a complex consisting of any two 1-non-intersecting cycles equals 0 too.  Our next goal is to show that the coefficient of the complex consisting of two simple cycles either 1-intersecting or not and a spike is 0.  Denote the complex with $\Delta$ and the two simple cycles with $\Delta_1$ and $\Delta_2$.  Let $\Delta'$ be the subcomplex derived form $\Delta$ by erasing the spike and let $\tilde{\Delta^1}$ be the set of complexes which is complement to $\Delta'$, i.e. such that $\Delta^1 = \Delta' \cup \tilde{\Delta^1}$.  Split formula (\ref{eq:qko}) as two sets of elements - one generated by all the subcomplexes obtained from $\Delta'$ and the second one - from $\tilde{\Delta^1}$.  If $\Delta$ consists of two 1-intersecting simple cycles and a spike, then $c_1(\Delta) =0$ follows from formula (\ref{eq:qko}) applied to $c_1(\Delta')$.  In case that we are dealing with two 1-non-intersecting simple cycles and a spike, then we can write formula (\ref{eq:qko}) as
$$c_1(\Delta) = 2 - c_1(\Delta') - c_1(\tilde{\Delta^1}) - c_1(\Delta^2) - \ldots - c_1(\Delta_1) -c_1(\Delta_2)$$
where all the coefficients except $c_1(\Delta_1)$ and $c_1(\Delta_2)$ are 0, so once again we have 
$$c_1(\Delta) =0.$$

Finally, we can extend the result for a complex of two simple cycles and k spikes by induction if we suppose that the coefficient of two cycles with $k-1$ spikes is zero, using exactly the same arguments as above, we get that $c_1(\Delta) = 0$.  To summarize everything known so far, we discovered that the coefficient of any complex of either 1-intersecting or not two simple cycles and $k>0$ of spikes is always 0.  By the same argument, we can take arbitrary complex with non-zero coefficient add a spike to it, use formula (\ref{eq:qko}), and get that the coefficient of the complex is 0.  Then, using induction, one can show that $c_1$ of the union of the non-zero coefficient complex and k spikes is again 0.

Our final aim is to expand the result for any number of simple cycles and spikes satisfying the previous restrictions, i.e. we will show that the coefficient of a complex consisting of any number of simple cycles, either 1-intersecting or not and any positive number of spikes is always 0, and also - the coefficient of a complex which is build up from two 1-non-intersecting complexes is also 0.  Then one can use the same argument which we use and extend the result from simple cycles to any complexes.  Once again, we get use of the simple but powerful method of induction, this time the induction will go on the number of component.  Suppose that the coefficients of any number of 1-non-intersecting cycles up to $p>1$ is zero and the coefficients of any number of 1-non-intersecting cycles up to some number $b_0^E$ with any number of spikes is also zero.  Then the proof is done in two steps.  First, take $\Delta$ to be $p+1$ 1-non-intersecting fundamental cycles.  Using formula (\ref{eq:qko}) we have that
$$c_1(\Delta) = p+1 - c_1(\Delta^1) - c_1(\Delta^2) - \ldots - \sum_{i=1}^{p+1} c_1(\Delta_i).$$
By the induction hypothesis, all $c_1(\Delta^i) = 0$, thus $c_1(\Delta)=0$.  Similarly as before, using induction again, we can show that $p+1$ cycles with any number of spikes have coefficient equal to 0.  

Thus we revealed everything about the structure of the polynomial which measures $b_0^E$.

\begin{thm}\label{thm:bettizero}
After reduction, the polynomial $b_0^E$ has the structure
$$b_0^E = p_1 + \ldots + p_n - \sum_{i<j} p_{ij} + \mbox{(terms of higher order)}$$
where the terms of higher order with nonzero coefficients are built by simple cycles having 1-intersections.  The coefficients of these complexes are given by formula (\ref{eq:qko}).  The coefficients of all other complexes are 0.
\end{thm}

There are two elementary operations which build up any complex.  The first one $\vartriangle$ is symmetric difference and the way it works is to take the union of two complexes and subtract the 1-dimensional cells from the common intersection.  The second one $\cup$ is the usual union.  We add subindex $1$ to the operations if we want to clarify that the operation applies at one 1-dimensional cells or $2$ if it applies at two 1-dimensional cells.

Shape of the figures that build up the complex does not matter as soon as they are topologically the same.

\begin{lemma}
For any complex $\Delta$ and a triangle $\Delta^t$,
$$c_1(\Delta \vartriangle_1 \Delta^t) = c_1(\Delta).$$
\end{lemma}
\textbf{Proof.}  Notice that $\vartriangle_1$ does not change $b_1$, i.e. $b_1(\Delta \vartriangle_1 \Delta^t) = b_1(\Delta)$ and let $b_1(\Delta)=n$.  Then using Theorem \ref{thm:govna} we have
\begin{equation}\label{eq:pederasi}
\sum_{\Delta^i \in \Delta \vartriangle_1 \Delta^t} c_1(\Delta^i) = \sum_{\Delta^i \in \Delta} c_1(\Delta^i).
\end{equation}
The proof uses induction on $b_1$.  If $b_1(\Delta)=1$, we already know that $c_1(\Delta \vartriangle_1 \Delta^t) = c_1(\Delta)$.  Suppose, it is true that $c_1(\Delta \vartriangle_1 \Delta^t) = c_1(\Delta)$ for all $\Delta$, s.t. $b_1(\Delta) < n$, then by equation (\ref{eq:pederasi}) and the fact that there is only one structure with $b_1 = n$, the proof follows.
\begin{flushright}
$\Box$
\end{flushright}

\begin{thm}
For any structure $\Delta$ and a triangle $\Delta^t$,
$$c_1(\Delta \cup_1 \Delta^t) = - c_1(\Delta).$$
\end{thm}
\textbf{Proof.}  We use a technique similar to Mayer-Vietoris sequence.  Decompose the complex $\Delta \cup_1 \Delta^t$ as a union of $\Delta$, $\Delta^t$ and $\Delta'$, where the last one is the set of complexes which are build by elementary operations on elements from both $\Delta$ and $\Delta^t$.  If we substitute the elements from  $\Delta$ and $\Delta^t$ in formula (\ref{eq:qko}), we get that 
\begin{equation}\label{eq:ako}
\sum_{\omega \in \Delta} c_1(\omega * \Delta^t) = 0
\end{equation}
where $*$ is any of the elementary operations.  Notice that if $\omega \in \Delta$ and $\Delta^t$ are 1-nonintersecting, then $c_1(\omega \cup_1 \Delta^t) = c_1(\omega \vartriangle_1 \Delta^t) = 0$.  The rest of the proof can be done by induction starting with $$c_1(\omega \cup_1 \Delta^t) = - c_1(\omega \vartriangle_1 \Delta^t) = - c_1(\omega) = -1,$$ when $\omega$ is a fundamental cycle.  Then using the induction hypothesis and formula (\ref{eq:ako}), it follows that $c_1(\Delta \cup_1 \Delta^t) = - c_1(\Delta).$
\begin{flushright}
$\Box$
\end{flushright}

\section{Term Structure of $b_n^E$}

Consider the following monomial
$$p_1 \ldots p_i p_{i+1} \ldots p_j p_{j+1} \ldots p_k p_{p_1 p_2} \ldots p_{p_{2m-1} p_{2m}} p_{p_{2m+1} p_{2m+2}} \ldots p_{p_{2n-1} p_{2n}} p_{q_1 q_2 q_3} \ldots p_{q_{3s-2} q_{3s-1} q_{3s}}$$
such that $m$ is the minimum number of edges which is covered by $s$ 2-faces, $i$ is the minimum number of points covered by these $m$ edges.  Also, there are another $n-m$ edges which cover $j-i$ points.  Thus, once again we are forced to choose the first $i$ together with the next $j-i$ points and another $m$ edges to exist.  So, we can choose only whether $k-j$ points, $n-m$ edges and $s$ faces exist or not.  This time we are interested in summing those monomials with coefficients - the first betti number of the monomial.  Then the following expression gives the coefficient of the above monomial
$$\sum_{t=0}^{k-j} (-1)^t \binom{k-j}{t} \sum_{\substack{| \nu | = 0 \\ \nu \in C(p_{2m+1}, \ldots, p_{2n})}}^{n-m} (-1)^{|\nu|} \sum_{\substack{| \omega | = 0 \\ \omega \in C(q_1, \ldots, q_{3s})}}^s (-1)^{|\omega|}-(k-t-n+| \nu | + s - |\omega| - b_0^{\omega} - b_2^{\omega})$$
where $t$ is the number of non-existing points, $|\nu|$ is the number of non-existing edges and $|\omega|$ is the number of non-existing faces and in the parenthesis you can discover Euler - Poincare formula for $b_1$.  Notice that $k-t-n+| \nu | - b_0^{\omega}$ does not depend on $\omega$ as $b_0^{\omega}$ does not depend on whether there is a face or not and instead we can write $b_0^{\nu}$.  So we can split the sum into three sums and evaluate each one of them
$$\begin{array}{lcr}
-\sum_{t=0}^{k-j} (-1)^t \binom{k-j}{t} \sum_{\substack{| \nu | = 0 \\ \nu \in C(p_{2m+1}, \ldots, p_{2n})}}^{n-m} (-1)^{|\nu|} ( \sum_{\substack{| \omega | = 0 \\ \omega \in C(q_1, \ldots, q_{3s})}}^s (-1)^{|\omega|}(k-t-n+| \nu | - b_0^{\nu}) & & \\
+  \sum_{\substack{| \omega | = 0 \\ \omega \in C(q_1, \ldots, q_{3s})}}^s (-1)^{|\omega|} (s-|\omega|) + \sum_{\substack{| \omega | = 0 \\ \omega \in C(q_1, \ldots, q_{3s})}}^s (-1)^{|\omega|+1}b_2^{\omega} )
\end{array}$$
Now, actually $k-t-n+| \nu | - b_0^{\nu}$ is exactly $-b_1^{\nu}$ and let's remove the combinatorial term wherever possible
$$\begin{array}{lcr}
-\sum_{t=0}^{k-j} (-1)^t \binom{k-j}{t} \sum_{\substack{| \nu | = 0 \\ \nu \in C(p_{2m+1}, \ldots, p_{2n})}}^{n-m} (-1)^{|\nu|} ( - b_1^{\nu} \sum_{| \omega | = 0 }^s (-1)^{|\omega|}\binom{s}{|\omega|} & & \\
+  \sum_{| \omega | = 0 }^s (-1)^{|\omega|} (s-|\omega|)\binom{s}{|\omega|} + \sum_{\substack{| \omega | = 0 \\ \omega \in C(q_1, \ldots, q_{3s})}}^s (-1)^{|\omega|+1}b_2^{\omega} ) & & \\
= -\sum_{t=0}^{k-j} (-1)^t \binom{k-j}{t} \sum_{\substack{| \nu | = 0 \\ \nu \in C(p_{2m+1}, \ldots, p_{2n})}}^{n-m} (-1)^{|\nu|} ( - b_1^{\nu} (1-1)^s & & \\
+ (-1)^{s-1} s \sum_{| \omega | = 0 }^s (-1)^{s-1-|\omega|} \binom{s-1}{|\omega|} + \sum_{\substack{| \omega | = 0 \\ \omega \in C(q_1, \ldots, q_{3s})}}^s (-1)^{|\omega|+1}b_2^{\omega} ) & & \\
= -\sum_{t=0}^{k-j} (-1)^t \binom{k-j}{t} \sum_{\substack{| \nu | = 0 \\ \nu \in C(p_{2m+1}, \ldots, p_{2n})}}^{n-m} (-1)^{|\nu|} ( - b_1^{\nu} (1-1)^s & & \\
+ (-1)^{s-1} s (1-1)^{s-1} + \sum_{\substack{| \omega | = 0 \\ \omega \in C(q_1, \ldots, q_{3s})}}^s (-1)^{|\omega|+1}b_2^{\omega} )
\end{array}$$
Once again we can evaluate the first sum as $b_1^{\nu}$ and $b_2^{\omega}$ do not depend on the number of non-existing points, so we get $k=j$ and so the coefficient is
$$-\sum_{\substack{| \nu | = 0 \\ \nu \in C(p_{2m+1}, \ldots, p_{2n})}}^{n-m} (-1)^{|\nu|} ( - b_1^{\nu} (1-1)^s + (-1)^{s-1} s (1-1)^{s-1} + \sum_{\substack{| \omega | = 0 \\ \omega \in C(q_1, \ldots, q_{3s})}}^s (-1)^{|\omega|+1}b_2^{\omega} ).$$
If $s=0$, then the coefficient is given by the first summand only, as all others are zero, i.e.
$$\sum_{\substack{| \nu | = 0 \\ \nu \in C(p_{2m+1}, \ldots, p_{2n})}}^{n-m} (-1)^{|\nu|} b_1^{\nu}$$
which is the coefficient $c_1$.  If $s>0$, we get the following coefficient
$$-\sum_{\substack{| \nu | = 0 \\ \nu \in C(p_{2m+1}, \ldots, p_{2n})}}^{n-m} (-1)^{|\nu|} ( (-1)^{s-1} s (1-1)^{s-1} + \sum_{\substack{| \omega | = 0 \\ \omega \in C(q_1, \ldots, q_{3s})}}^s (-1)^{|\omega|+1}b_2^{\omega} ).$$
We can evaluate the first sum and get that $n=m$ and the coefficient becomes
$$(-1)^{s} s (1-1)^{s-1} + \sum_{\substack{| \omega | = 0 \\ \omega \in C(q_1, \ldots, q_{3s})}}^s (-1)^{|\omega|}b_2^{\omega}.$$
When $s=1$ the second sum is zero and the coefficient is -1.  Similarly to the discussion above, define the coefficient $c_2$ to be
$$c_2(\Delta) = \sum_{\substack{| \omega | = 0 \\ \omega \in C(q_1, \ldots, q_{3s})}}^s (-1)^{|\omega|}b_2^{\omega}$$
if $b_2(\Delta)$ is non-zero.

In general, we can apply exactly the same procedure for calculating $b_{k-1}^E$ and get that each k-cell comes with sign equal to $(-1)^{k-1}$ and for $b_k(\Delta) \neq 0$ we have
$$c_k(\Delta) = \sum_{\substack{| \omega | = 0 \\ \omega \in C(q_1, \ldots, q_{ks})}}^{s} (-1)^{|\omega|}b_k^{\omega}.$$

Similarly to the results in the lower case, we can prove
\begin{cons} 
\begin{equation}\label{eq:kuku}
b_k(\Delta) = \sum_{i=0}^n c_k(\Delta^i)
\end{equation}
and so there are no terms corresponding to $m$-spike and $m$-nonintersecting cycles, where $m \leq k$.
\end{cons}

\begin{thm}\label{thm:bettin}
The polynomial $b_n^E$ has the following structure
$$b_n^E = p_n^E + (-1)^{n+1} d_{n+1} + p_{n+1}^E$$
where $p_i^E$ is the polynomial corresponding to the higher order terms of i-cell and $d_{n+1}$ is the polynomial corresponding to the (n+1)-cells.
\end{thm}

\section{Explicit Definitions}

\begin{defin} 
Let the random variable $\mathbb{C}_n$ measures the number of n-cells in the complex $X$, i.e. $\mathbb{C}_n = e_1 + \dots + e_n$, $e_i \in \{ 0, 1 \}$ where each $e_i$ is itself a random variable which corresponds to a real n-cell in the complex of all possible configurations of subcomplexes.  Define the expected number of n-cells $\mathcal{C}_n^E$ to be
$$\mathcal{C}_n^E = \sum c_j P(c_j) = E(\mathbb{C}_n)$$
where $P(c_j) = P(\mathbb{C}_n = c_j)$ and $c_j$ is the number of n-cells.
\end{defin}

Similarly, one can define the rank of the expected n-cycles and n-boundaries.

\begin{defin} 
Let the random variable $\mathbb{Z}_n$ measures the number of n-cycles in the complex $X$, i.e. $\mathbb{Z}_n = k_1 + \dots + k_m$, $k_i \in \{ 0, 1 \}$ where each $k_i$ is itself a random variable which corresponds to a real n-cycle in the complex.  Define the expected number of n-cycles $\mathcal{Z}_n^E$ to be
$$\mathcal{Z}_n^E = \sum z_j P(z_j)= E(\mathbb{Z}_n)$$
where $c_j$ is the number of n-cycles.
\end{defin}

\begin{defin} 
Let the random variable $\mathbb{B}_n$ measures the number of n-boundaries in the complex $X$, i.e. $\mathbb{B}_n = l_1 + \dots + l_p$, $l_i \in \{ 0, 1 \}$ where each $l_i$ is itself a random variable which corresponds to a real n-boundary.  Define the expected number of n-boundaries $\mathcal{B}_n^E$ to be
$$\mathcal{B}_n^E = \sum b_j P(b_j)= E(\mathbb{B}_n)$$
where $b_j$ varies over the number of n-boundaries.
\end{defin}

From probability theory, we know that
\begin{equation}\label{eq:imp}
E(X+Y) = E(X) + E(Y)
\end{equation}
thus it directly follows that
$$\mathcal{C}_n^E = \mathcal{Z}_n^E + \mathcal{B}_{n-1}^E.$$

Whatever $H_n^E = \frac{Z_n^E}{B_n^E}$ is supposed to be, define the n-th betti number of a random complex to be
$$b_n^E = \mathcal{Z}_n^E - \mathcal{B}_n^E.$$

\begin{defin}
Define expected Euler characteristic to be
$$\chi^E = \sum_j (-1)^j \mathcal{C}_j^E.$$
\end{defin}

Applying equation \eqref{eq:imp} to the last definition and taking in mind the cancellations which occur by the minus sign in each consecutive summand in the above equation, we get

\begin{thm}\label{thm:euler}Expected Euler-Poincare Formula.
$$\chi^E = \sum_j (-1)^j b_j^E.$$
\end{thm}

The last theorem can be easily verified by Theorems \ref{thm:bettizero} and \ref{thm:bettin}.

From Theorem \ref{thm:euler} it follows

\begin{thm}  $$\chi^E (A \cup B) = \chi^E (A) + \chi^E (B) - \chi^E (A \cap B)$$
for random complexes $A$ and $B$.
\end{thm}

In general, one can build a complete stochastic homology theory by redefining the classical homology axioms in stochastic homology sense.  The classical homology theory sits in the stochastic and one can get most of the classical results valid for the stochastic homology as well.  However, the most important result for calculation, e.g. Mayer-Vietoris exact sequence, does not hold.  The reason for that is that the polynomial of higher terms $p_n^E$ of union of complexes contains terms from both subcomplexes, and thus knowing the expected betti numbers of the two subcomplexes is simply not enough to generate the expected betti number of the union.

\section{Computation, Algorithm and Experimental Results}

A quick look at Theorems \ref{thm:bettizero} and \ref{thm:bettin} reveals that calculation of expected betti numbers in general seems to be impossible for big set of points as the number of summands in $b_n^E$ grows rapidly.  However, one can notice that each polynomial representing $b_n^E$ has several symmetries in it.  It turns out that those polynomials are stabilized by a subgroup $S_m \hookrightarrow S_{\binom{m}{n}}$, where $m$ is the number of 0-dimensional cells of the complex.  This is not surprising as any shift of two vertices generate a basis element for $S_m$.  Such a polynomial is called to be almost symmetric and my hope was to try to represent it as product of linear terms over the field of complex numbers.  Unfortunately, my attempts did not give a positive result.

Thus, the only reasonable calculation of this polynomials can occur if we fix all probabilities of cells of dimension greater than 1 to be equal.  Let's concentrate on the calculation of the polynomial $p_1^E$, i.e. $p_{12} = \ldots = p_{(m+1) m} = x$.  For $m=4$, we have $$p_1^E = 2 x^6 - 6 x^5 + 3 x^4 + 4 x^3$$
and for $m=5$
$$p_1^E = -6 x^{10} + 40 x^9 - 105 x^8 + 130 x^7 - 60 x^6 - 18 x^5 + 15 x^4 + 10 x^3.$$
In general, the polynomial $p_n^E$ is of degree $\binom{m}{n+1}$ and this polynomial is easily computable at least for small $n$ even when $m$ is large.  Thus one can calculate expected betti numbers for a Vietoris-Rips complex.

The next question is how to assign probabilities to the cells when you have no statistics available for their distribution.  One geometric way is to look at the distance between the center of mass of the cell and the closest point to it available from the data.  Call this distance $r_d$ and denote with $r_m$ the distance between the center of mass and a vertex of the cell.  There are many ways to assign probability using this data.  I was looking at the circles centered at the center of mass with radii $r_d$ and $r_m$, thus the first probability which I constructed was
$$p = 1 - \frac{\pi r_d^2}{ \pi r_m^2} = 1 - \frac{r_d^2}{r_m^2}.$$
It turns out that this probability shrinks the gaps between the points because the volume of a circle is contained at the exterior and so it was not what I was looking for.  However, it gives a good idea for a better probability
$$p = 1 - \left(\frac{r_d}{r_m}\right)^{1/2}$$
which works well for random complexes.  The whole class
$$p = 1 - \left(\frac{r_d}{r_m}\right)^{1/k}, \quad p = 1 - \left(\frac{r_d}{r_m}\right)^{k}, \quad k \in \NN$$
seems to be helpful for different types of data.

Once the probabilities are assigned, we generate monomials corresponding to n-cells of different degrees.  Each monomial has a unique coefficient $c_n$ and is distributed in the random complex by the action of the stabilizer subgroup.  Thus, it is only needed to calculate the order $o_n$ of the orbit group, which can be done by combinatorial methods.  The coefficients $c_n$ can be either zero or non-zero.  

The first time I discovered the magnificent properties of these numbers, I spontaneously called them "magical" coefficients.  Though I could't prove anything more than the result cited above, the following experimental results reappear for different complexes.  The first one is to decompose combinatorically a full cycle except one edge.  The sum of all coefficients is zero.

\begin{center}
\includegraphics[scale=0.8]{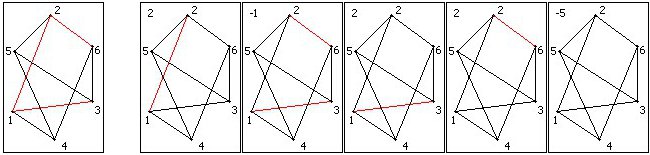}
\begin{center}
\small{Decomposition of cycle}
\end{center}
\end{center}

We decompose the cycle 1-3-6-2, but leave the edge 3-6 intact.  The slide on the left shows in red which edges we decompose.  The red edges in the right slides denote missing edges.  The coefficient is printed in the left upper corner. The slides with zero coefficients are not shown above. 

The next result which reappear again and again is if you decompose combinatorically the edges at a vertex.  The sum of all such coefficients is also zero.

\vspace{10mm}

\begin{center}
\includegraphics[scale=0.8]{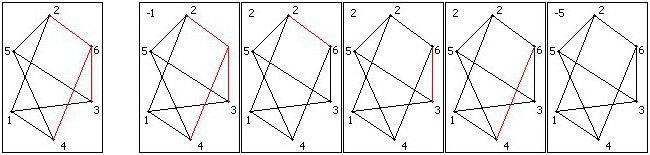}
\begin{center}
\small{Decomposition of edges at a vertex}
\end{center}
\end{center}

Here we decompose the edges at vertex 6 combinatorically.

The reason to call these coefficients "magical" is that you can get even more interesting results.  You can add more edges to the above decompositions and still get that the sum of the coefficients is zero.

\begin{center}
\includegraphics[scale=0.8]{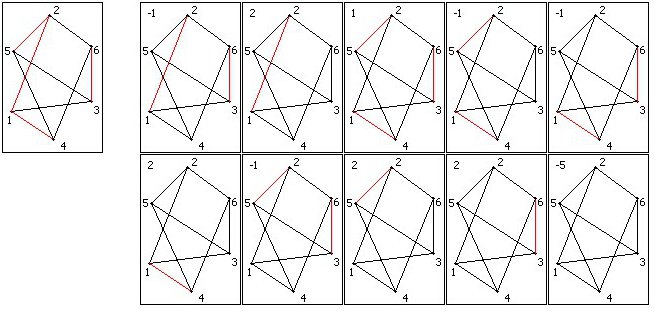}
\begin{center}
\small{Decomposition of a cycle and an edge}
\end{center}
\end{center}

And we can continue adding edges to be decomposed, and each time the sum of the coefficients appear to be zero.

\begin{center}
\includegraphics[scale=0.8]{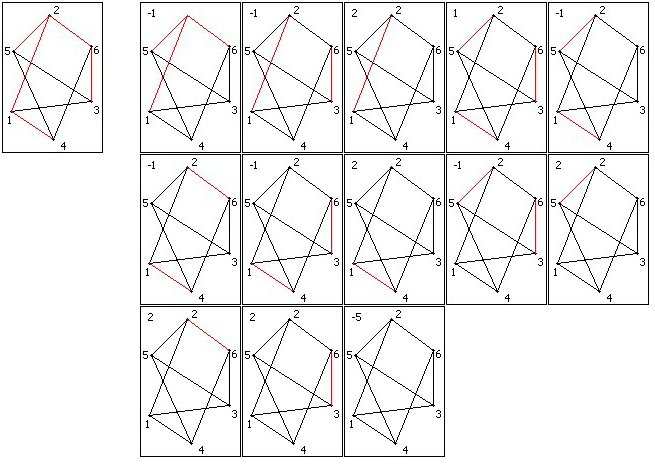}
\begin{center}
\small{Decomposition of cycle and two edges}
\end{center}
\end{center}

In fact, these coefficients are so designed that if we set all probabilities of cells equal to one, then no matter what complex we have the top coefficient, e.g. the one corresponding to the maximal monomial in this case, always wipes out all other coefficients and leaves the betti number only.  This is the reason why Mayer-Vietoris holds in the classical case, but not in the random.

Back to the discussion of the algorithm - the zero coefficients are associated with either a complex which has a spike, or a complex which is built up from two or more cycles which have zero n-intersection.  The spikes can be detected by counting the degrees of the vertices.  We say that a vertex has degree $h$ if h cells cover that vertex.  Thus there is a $n$-spike in the complex if there is a vertex with degree $n$.  The $n$-intersection function is a little bit more complicated.  We remove each $n$-cycle which has $n$-intersection with another cycle.  At the end of the procedure, if there is a vertex of degree greater than $n+1$, then the complex does not satisfy the $n$-intersection test.  If the monomial does not fail both tests above, i.e. $c_n$ is nonzero, then it can be computed using formula (\ref{eq:kuku}) recursively.  The cheapest way to do it is perhaps using binary tree with a union cell.  However, I found it also efficient using binary tree structure with one coefficient cell.  To each boolean word of zeros and ones as 101010110011, we can assign direction in the binary tree, left for 0 and right for 1, and it is very easy to work with the coefficients this way.  There are two more conditions which speed up the process by looking at the reduced form, which we get by erasing a vertex of degree $n+1$ who has n+1 neighbors of degree $n+1$ and adding new $n$-cell between its $n+1$ neighbors.  So if the reduced form is different from the initial form of the monomial, then we can simply copy the coefficient from the database.  Finally, if the monomial has a vertex of degree $n$, then we can erase the $n$-cycle passing through that vertex, adjust the remaining degrees and copy the coefficient of the new monomial from the database with negative sign.   Now we can write $k_i = \sum o_n c_n$ and the formula for $p_n^E$ can be written as
$$p_n^E = \sum_i k_i x^i.$$
Once we have the formula for $p_n^E$ calculated once and forever, we can simply evaluate it using nested sequence of multiplications, that is
$$p_n^E = x^3 ( \ldots ((k_{\binom{m}{n}} x + k_{\binom{m}{n}-1})x + k_{\binom{m}{n}-2})x + \ldots + k_3).$$
\begin{center}
\begin{tabular}{|l|}
\hline
for(i=Binom\{m,n\}; i$\geq$3; i- -)\\
\quad do\\
\quad \quad monomial = GenerateNewMonomial(i);\\
\quad \quad \quad if(Spike(monomial)==false)\\
\quad \quad \quad \quad if(NIntersection(ContractMonomial(monomial))==false)\\
\quad \quad \quad \quad \quad newMonomial = ReduceMonomial(monomial);\\
\quad \quad \quad \quad	\quad oN = CalculateOrbitOrder(newMonomial);\\
\quad \quad \quad \quad \quad if(newMonomial != monomial)\\
\quad \quad \quad \quad \quad \quad cN = ReadCoefficient(newMonomial);\\
\quad \quad \quad \quad \quad if(vertexOfDegreeNExists(monomial)==true)\\
\quad \quad \quad \quad \quad \quad cN = - ReadCoefficient(EraseCycle(newMonomial));\\
\quad \quad \quad \quad \quad else\\
\quad \quad \quad \quad \quad \quad cN = CalculateCoefficientUsingFormula(monomial);\\
\quad \quad \quad \quad \quad kN+=oN*cN;\\
\quad while(monomial!=NULL);\\
\quad pN+=(pN+kN)*x;\\
\quad kN = 0;\\
return pN*x*x*x;\\
\hline
\end{tabular}
\begin{center}
\small{Algorithm for computing expected n-th betti number}
\end{center}
\end{center}

\section{Applications.  Coverage Problems and Large High-Dimensional Data}

There are perhaps thousands of problems where the model of stochastic homology can be applied.  We will focus on two of them - coverage problems and high-dimensional data with noise.

Coverage problem is set when a region is given together with subsets in the region and the question is whether these sets cover the entire region.  Perhaps the cheapest way to solve such problem is using algebraic topology methods and restate the question of whether the compactified union of subsets has the same homology as the compactified region.  As an example, you can take the covering sets to be circles and associate the centers of the circles with cell phone towers, and the circles themselves with the regions where the signal is clear.  

\begin{center}
\includegraphics[scale=0.32]{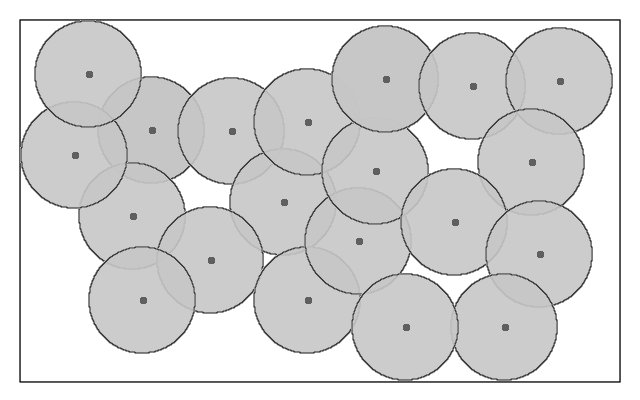}
\begin{center}
\small{Coverage Problem}
\end{center}
\end{center}

The zeroth betti number of the union of circles gives whether the region is connected, i.e. whether each cell tower can be connected to any other via a path.  The first betti number gives the number of cycles which cannot be contracted, i.e. the holes in the coverage.  Nowadays applied algebraic topology methods, such as persistence homology can answer such problems.

In real life, however, thing are not that simple.  Noise usually appears on the map of cell tower coverage.  Bad weather, interception and component failure are the most common factors and the maps in real time looks more like an animation of appearing and disappearing contracting and extracting regions which might not even be circles.  In this complicated formulation, we only have some statistics of the work times of the towers and of the quality of the signal between towers for some period of time.  This is enough to construct random complex and apply stochastic homology methods in order to find some estimate of whether the network had good signal coverage over time.

Another problem which can be associated with stochastic homology is the homology type of high-dimensional data.  You can think of such data as the set formed by a large collection of pictures, where each picture is mapped to a point in a very high dimensional space $\RR^d$, $d$ being the number of pixels in the picture.  Other example is financial data - each point is associated with a vector with components the price of each product in the market.  The evolution in time creates a path, which might be considered to lie on a manifold or a set with variance in the normal bundle.  We can approximate the manifold locally with the best fit hyperplane, or the one which minimizes the square distances.   In both case, usual persistent homology methods can be used, but the use of Voronoi cells, weak witnesses, etc. in order to create a complex restricts us to a certain level of certainty.  The use of stochastic homology breaks these limits and allows us to use a level of uncertainty in our model.  We can either guess that a point of the manifold exists, even though there is no evidence in the data that it is true, or we can simply thin the manifold in the regions where it is not dense.

\end{document}